\documentclass[12pt]{article}
\usepackage{amsfonts}
\begin{document}
\newsymbol\rtimes 226F
\newsymbol\ltimes 226E
\newcommand{\text}[1]{\mbox{{\rm #1}}}
\newcommand{\Rep}{\text{Rep}}
\newcommand{\End}{\text{End}}
\newcommand{\gd}{\delta}
\newcommand{\itms}[1]{\item[[#1]]}
\newcommand{\nin}{\in\!\!\!\!\!/}
\newcommand{\g}{{\bf g}}
\newcommand{\sub}{\subset}
\newcommand{\cntd}{\subseteq}
\newcommand{\go}{\omega}
\newcommand{\Pa}{P_{a^\nu,1}(U)}
\newcommand{\fx}{f(x)}
\newcommand{\fy}{f(y)}
\newcommand{\gD}{\Delta}
\newcommand{\gl}{\lambda}
\newcommand{\gL}{\Lambda}
\newcommand{\half}{\frac{1}{2}}
\newcommand{\sto}[1]{#1^{(1)}}
\newcommand{\stt}[1]{#1^{(2)}}
\newcommand{\Z}{\hbox{\sf Z\kern-0.720em\hbox{ Z}}}
\newcommand{\singcolb}[2]{\left(\begin{array}{c}#1\\#2
\end{array}\right)}
\newcommand{\ga}{\alpha}
\newcommand{\gb}{\beta}
\newcommand{\gga}{\gamma}
\newcommand{\ul}{\underline}
\newcommand{\ol}{\overline}
\newcommand{\qed}{\kern 5pt\vrule height8pt width6.5pt depth2pt}
\newcommand{\Lrraro}{\Longrightarrow}
\newcommand{\Nb}{|\!\!/}
\newcommand{\NN}{{\rm I\!N}}
\newcommand{\bsl}{\backslash}
\newcommand{\gt}{\theta}
\newcommand{\op}{\oplus}
\newcommand{\C}{{\bf C}}
\newcommand{\Q}{{\bf Q}}
\newcommand{\Op}{\bigoplus}
\newcommand{\CR}{{\cal R}}
\newcommand{\tr}{\bigtriangleup}
\newcommand{\grr}{\omega_1}
\newcommand{\ben}{\begin{enumerate}}
\newcommand{\een}{\end{enumerate}}
\newcommand{\ndiv}{\not\mid}
\newcommand{\bab}{\bowtie}
\newcommand{\hal}{\leftharpoonup}
\newcommand{\har}{\rightharpoonup}
\newcommand{\ot}{\otimes}
\newcommand{\OT}{\bigotimes}
\newcommand{\bwe}{\bigwedge}
\newcommand{\gep}{\varepsilon}
\newcommand{\gs}{\sigma}
\newcommand{\rbraces}[1]{\left( #1 \right)}
\newcommand{\bbox}{$\;\;\rule{2mm}{2mm}$}
\newcommand{\sbraces}[1]{\left[ #1 \right]}
\newcommand{\bbraces}[1]{\left\{ #1 \right\}}
\newcommand{\OO}{_{(1)}}
\newcommand{\T}{{\mathcal T}}
\newcommand{\FF}{_{(3)}}
\newcommand{\minus}{^{-1}}
\newcommand{\CV}{\cal V}
\newcommand{\CVs}{\cal{V}_s}
\newcommand{\un}{U_q(sl_n)'}
\newcommand{\on}{O_q(SL_n)'}
\newcommand{\slq}{U_q(sl_2)}
\newcommand{\olq}{O_q(SL_2)}
\newcommand{\UU}{U_{(N,\nu,\go)}}
\newcommand{\HH}{H_{n,q,N,\nu}}
\newcommand{\ct}{\centerline}
\newcommand{\bs}{\bigskip}
\newcommand{\qua}{\rm quasitriangular}
\newcommand{\ms}{\medskip}
\newcommand{\noin}{\noindent}
\newcommand{\mat}[1]{$\;{#1}\;$}
\newcommand{\raro}{\rightarrow}
\newcommand{\map}[3]{{#1}\::\:{#2}\raro{#3}}
\newcommand{\alg}{{\rm Alg}}
\def\newtheorems{\newtheorem{theorem}{Theorem}[subsection]
                 \newtheorem{cor}[theorem]{Corollary}
                 \newtheorem{prop}[theorem]{Proposition}
                 \newtheorem{lemma}[theorem]{Lemma}
                 \newtheorem{defn}[theorem]{Definition}
                 \newtheorem{Theorem}{Theorem}[section]
                 \newtheorem{Corollary}[Theorem]{Corollary}
                 \newtheorem{Proposition}[Theorem]{Proposition}
                 \newtheorem{Lemma}[Theorem]{Lemma}
                 \newtheorem{Definition}[Theorem]{Definition}
                 \newtheorem{Example}[Theorem]{Example}
                 \newtheorem{Remark}[Theorem]{Remark}
                 \newtheorem{claim}[theorem]{Claim}
                 \newtheorem{sublemma}[theorem]{Sublemma}
                 \newtheorem{example}[theorem]{Example}
                 \newtheorem{remark}[theorem]{Remark}
                 \newtheorem{question}[theorem]{Question}
                 \newtheorem{Question}[Theorem]{Question}
                 \newtheorem{conjecture}{Conjecture}[subsection]}
\newtheorems
\newcommand{\Proof}{\par\noindent{\bf Proof:}\quad}
\newcommand{\dmatr}[2]{\left(\begin{array}{c}{#1}\\
                            {#2}\end{array}\right)}
\newcommand{\doubcolb}[4]{\left(\begin{array}{cc}#1&#2\\
#3&#4\end{array}\right)}
\newcommand{\qmatrl}[4]{\left(\begin{array}{ll}{#1}&{#2}\\
                            {#3}&{#4}\end{array}\right)}
\newcommand{\qmatrc}[4]{\left(\begin{array}{cc}{#1}&{#2}\\
                            {#3}&{#4}\end{array}\right)}
\newcommand{\qmatrr}[4]{\left(\begin{array}{rr}{#1}&{#2}\\
                            {#3}&{#4}\end{array}\right)}
\newcommand{\smatr}[2]{\left(\begin{array}{c}{#1}\\
                            \vdots\\{#2}\end{array}\right)}

\newcommand{\ddet}[2]{\left[\begin{array}{c}{#1}\\
                           {#2}\end{array}\right]}
\newcommand{\qdetl}[4]{\left[\begin{array}{ll}{#1}&{#2}\\
                           {#3}&{#4}\end{array}\right]}
\newcommand{\qdetc}[4]{\left[\begin{array}{cc}{#1}&{#2}\\
                           {#3}&{#4}\end{array}\right]}
\newcommand{\qdetr}[4]{\left[\begin{array}{rr}{#1}&{#2}\\
                           {#3}&{#4}\end{array}\right]}

\newcommand{\qbracl}[4]{\left\{\begin{array}{ll}{#1}&{#2}\\
                           {#3}&{#4}\end{array}\right.}
\newcommand{\qbracr}[4]{\left.\begin{array}{ll}{#1}&{#2}\\
                           {#3}&{#4}\end{array}\right\}}

\title{Geometric crystals and set-theoretical solutions to the quantum
Yang-Baxter equation} 
\author{Pavel Etingof\\
Massachusetts Institute of Technology\\
Department of Mathematics, Rm 2-176\\
Cambridge, MA 02139\\
{\rm email: etingof@math.mit.edu}
}
\maketitle

\section{Introduction}

The notion of a geometric crystal was introduced and developed 
recently in \cite{BeK}, motivated by the needs of representation 
theory of p-adic groups. It is shown in \cite{BrK,BeK} that some particular
geometric crystals give rise to interesting set-theoretical solutions 
$R$ of the quantum Yang-Baxter equation, which satisfy the 
involutivity (or unitarity) condition $R^{21}R=1$
(more precisely, they give rise to ``rational'' set theoretical solutions, 
see below).

On the other hand, involutive set-theoretical solutions 
of the quantum Yang-Baxter equation are studied in \cite{ESS}, 
following a suggestion of Drinfeld, \cite{Dr}. It turned out in \cite{ESS} 
that there is an especially nice theory of such solutions if an additional 
nondegeneracy condition is satisfied. In particular, 
in this case one can define the structure group $G_R$ 
of a solution $(X,R)$, which acts on $X$. 
The image $G_R^0$ of $G_R$ in $\text{Aut}(X)$ is called the reduced 
structure group. The complexity of this group, in a sense, 
characterizes the complexity of the solution $(X,R)$. 

In this note we show that the maps $(X,R)$ arising from the geometric 
crystals of \cite{BrK},\cite{BeK} are nondegenerate, and 
give a new proof that they satisfy the quantum Yang-Baxter equation 
and the involutivity condition. Then we calculate 
the reduced structure group of $(X,R)$ and show that it is a 
subgroup of the group $PGL_n(\C(\lambda))$. 
We also give a new, direct proof of Theorem 8.9 of \cite{BrK}. 

{\bf Remark.} There is a seemingly technical, but in fact fundamental 
point that needs to be stressed. The nondegenerate solutions 
coming from geometric crystals live {\it not} in the category of sets 
(like the solutions from \cite{ESS}) but in the category of 
irreducible algebraic varieties, where morphisms are {\it birational} 
maps. In particular, they {\it cannot} be viewed as 
usual nondegenerate set-theoretical solutions, studied in \cite{ESS};
they are similar but yet essentially new, more complicated objects. 
As a result, the theory of \cite{ESS} has to be generalized 
to the ``rational'' case, to be applicable to geometric crystals. 
This generalization is not straightforward, and it is not entirely 
clear which parts of the theory survive and which do not. 
This seems to be an interesting problem for future research. 

{\bf Acknowledgments.} The author thanks A.Braverman and D.Kazhdan 
for useful discussions. This work was partially supported by the 
NSF grant DMS-9988796, and conducted in part for the Clay Mathematics 
Institute. 

\section{Rational set-theoretical R-matrices}

Let us generalize the setting of \cite{ESS} to the 
``rational'' case. 

Let $X$ be an irreducible algebraic variety over $\C$. 
Let $R$ be a birational isomorphism of $X\times X$ to itself. 

Let us write $R(x,y)$, $x,y\in X$, as 
\begin{equation}
R(x,y)=(f_y(x),g_x(y)).
\end{equation}
The rational mappings $f_z,g_z:X\to X$ are well defined for generic 
$z\in X$. 

\begin{Definition}  $R$ 
is said to be nondegenerate if $f_z,g_z$ are 
defined and are birational isomorphisms 
for all $z$. 
\end{Definition}

For example, $R=1$ is nondegenerate but $R=P$ (the permutation) is not. 

\begin{Definition} $R$ 
is said to be involutive if $R^{21}R=1$. 
\end{Definition}

The quantum Yang-Baxter equation for $R$ is the equation 
\begin{equation}
R_{12}R_{13}R_{23}=R_{23}R_{13}R_{12}
\end{equation}
between birational isomorphisms of $X^3$. 

\begin{Definition} A solution $R$ of the quantum Yang-Baxter equation 
will be called a rational set-theoretical R-matrix. 
\end{Definition}

For an irreducible variety $Y$, let $\text{Bir}(Y)$ be the group of birational 
automorphisms of $Y$. 
An involutive $R$-matrix $R$ defines, for each $N$, a homomorphism 
$\rho_N^R:S_N\to \text{Bir}(X^N)$, given by $\rho((i,i+1))=
PR_{i,i+1}$. For example, for $R=1$ this is the usual 
action by permutations, and for $R=P$ the trivial action. 

\begin{Proposition}
If $R$ is nondgenerate then $\rho_N^R$ is conjugate 
to the usual action $\rho_N^1$. Namely, 
$\rho_N^R=J_N^{-1}\rho_N^1J_N$, where 
$J_N\in \text{Bir}(X^N)$ is given by 
\begin{equation}
J_N(x_1,...,x_N)=(f_{x_N}...f_{x_2}(x_1),...,f_{x_N}(x_{N-1}),x_N).
\end{equation}
\end{Proposition} 

\begin{Proof} As in \cite{ESS}, Section 1. 
\end{Proof} 

Let $R$ be an involutive nondgenerate rational set-theoretical R-matrix. 
Let $U\subset X^2$ be the domain of definition of $R$, i.e. the largest 
open set such that $R$ is regular on $U$. Let $G_R$ be the group 
generated by the points of $X$ with defining relations
 \begin{equation}
xy=y'x'
\end{equation} 
if $(x,y)\in U$ and $R(x,y)=(x',y')$. The group $G_R$ is called 
the structure group of $(X,R)$. 

\begin{Proposition} The group $G_R$ acts on $X$ by birational 
transformations in two ways: 
$z\to f_z^{-1}$ and $z\to g_z$.  
\end{Proposition}

\begin{Proof} As in \cite{ESS}, Section 2. 
\end{Proof} 

Let us denote the image in $\text{Bir}(X)$ of the first action by 
$G_R^+$ and of the second one by $G_R^-$. These groups will 
be called reduced structure groups of $(X,R)$. 

\section{Rational set-theoretical R-matrices 
arising from geometric crystals.}

Let us recall the construction of \cite{BrK,BeK} (see \cite{BrK}, 
Section 8). 

Let $\T=(\C^*)^n$. Elements of 
$\T$ will be written as $\bold t=[t_1,...,t_n]$, $t_i\in \C^*$. 

For $k\ge 2$ let
\begin{equation}
\Delta_k(\bold x,\bold y)=x_1...x_{k-1}+x_1...x_{k-2}y_k+...
+y_2...y_k.
\end{equation}
(we agree that $\Delta_1=1$). 
Let 
\begin{equation}
\eta(\bold x,\bold y)=\frac{x_1...x_n-y_1...y_n}{\Delta_n(\bold x,\bold y)}
\end{equation}

The following proposition is easy. 

\begin{Proposition}\cite{BrK,BeK}
 There exists a unique rational map
$R:\T^2\to \T^2$,
\begin{equation}
R(\bold x,\bold y)=
(\bold x',\bold y'), 
\end{equation}
such that
\begin{equation}
x_1'...x_k'=y_1...y_k+\Delta_k(\bold x,\bold y)\eta(\bold x,\bold y),
\end{equation}
and 
\begin{equation}
x_i'y_i'=x_iy_i, 
\end{equation}
\begin{equation}
\prod_{i=1}^n x_i=
\prod_{i=1}^n x_i'.
\end{equation}
\end{Proposition}

Our first result is 

\begin{Proposition}\label{res1} $R$ is nondegenerate. 
\end{Proposition}

This proposition is proved in the next section. 

Moreover, one has the following result, which is 
a part of Theorem 8.9 of \cite{BrK}, proved in \cite{BeK}
using geometric crystals. 

\begin{Theorem} \label{Rmat}
$R$ is an involutive rational set-theoretical R-matrix.  
\end{Theorem}

We give a direct proof of this theorem in Section 5. 
We note that such a proof was known to the authors of \cite{BrK}, but 
as far as we know, it is unpublished. 

\begin{Corollary} The $S_N$ action on $\T^N$ defined by $R$
is conjugate to the usual one.
\end{Corollary}

{\bf Remark.} Since in the construction of $R$ one does not use subtraction
(see \cite{BeK}), one can regard $R$ as a usual (i.e. not ``rational'') 
set-theoretical R-matrix defined 
on the set of points of $\T$ with positive real 
coordinates. Another (essentially, equivalent) way to turn R into a usual 
set-theoretical R-matrix is to use tropicalization
(see \cite{BeK} and references therein). However, 
it is easy to see that these usual R-matrices will not be 
nondegenerate. We don't expect that it is possible 
to interpret R as a usual (not ``rational'')
{\it nondegenerate} set-theoretical R-matrix, 
and believe that it is essential to use the ``rational'' generalization 
to attain nondegeneracy. 
\vskip .05in

Our second result is 

\begin{Proposition}\label{res2} The groups $G_R^+$, $G_R^-$ are 
subgroups of $PGL_n(\C(\lambda))$. 
These subgroups are isomorphic to each other. 
\end{Proposition} 

This proposition is proved in Section 5. 
In the proof, we actually describe 
the groups $G_R^+,G_R^-$ explicitly. 

{\bf Remark.} We want to stress that we don't understand 
the meaning of the groups $G_R^\pm$ and of Proposition \ref{res2}
in the context of geometric crystals. We also don't know
how to generalize the main construction of \cite{ESS}
(the bijective cocycle construction) to the rational case, so that it would produce the solutions considered here. 
It seems that such a generalization would be interesting and useful. 
\vskip .05in

Finally, we give a direct proof of Theorem 8.9 of \cite{BrK}
(according to \cite{BrK}, a proof of this result 
without the theory of geometric crystals was previously unavailable).
More specifically, we prove the following statement, 
which is what Theorem 8.9 of \cite{BrK} claims 
in addition to Theorem \ref{Rmat}. 

Let $M=(\C^*)^{mn}$ be the set of $n$ by $m$ matrices with entries from 
$\C^*$. Using the R-matrix $R$ defined in Section 3, 
we can define two symmetric group actions on $M$: one of 
$S_m$ (on columns) and another of $S_n$ (on rows). 

\begin{Theorem} \label{comm}
The actions of $S_n$ and $S_m$ on $M$ commute with each other. 
\end{Theorem}

This theorem is proved in Section 6. 

\section{Proof of Proposition \ref{res1}}

Let us introduce functions on $T_i: \T\to \C^*$ by 
$T_0(\bold t)=1$, and $T_{i}=T_{i-1}t_{i\text{ mod }n}$, where 
$i\text{ mod }n$ takes values $1,...,n$.  
In this way, $T_j(\bold t)$ is defined for all integer $j$.

We have
$\T=\cup_{c\in \C^*}\T_c$, where $\T_c=\lbrace{\bold t\in \T|
T_n=c\rbrace}$. It is clear that 
for (almost) all $c$, the map $f_{\bold z}$,
$g_{\bold z}$ map $\T_c$ to $\T_c$.

Define two compactifications of $\T_c$, the projective spaces 
$P_+$ and $P_-$. To do this, it is sufficent to define two open embeddings 
$j_+,j_-:\T_c\to \C P^{n-1}$. We define them by
\begin{equation}
j_+(\bold t)=(1,T_1(\bold t),...,T_{n-1}(\bold t));\
j_-(\bold t)=(T_1(\bold t)^{-1},...,T_{n-1}(\bold t)^{-1},T_n(\bold t)^{-1}=
c^{-1});\
\end{equation}
The projective coordinates on 
$P_+,P_-$ will be denoted by 
$\bold Z=(Z_1,...,Z_n)$. 

Now let us take $\bold x\in \T_a$, $\bold y\in \T_b$, and 
calculate $f_{\bold y}:\T_a\to \T_a$, $g_{\bold x}: \T_b\to \T_b$. 
Let $X_i=T_i(\bold x),Y_i=T_i(\bold y)$. 

\begin{Proposition} 
(i) The map $f_{\bold y}: \T_a\to \T_a$ 
extends to a projective transformation of 
$P_+$ given by
\begin{equation}
\bold Z\to A_f(\bold y,a)\bold Z,
\end{equation} 
where $\bold Z$ is understood as a column vector, and 
$A_f$ is the matrix given by 
\begin{equation}
A_f(\bold y,a)_{ij}=Y_{i-1}Y_{j}^{-1}a, i>j;\  
A_f(\bold y,a)_{ij}=Y_{i-1}Y_{j}^{-1}b, i\le j,  
\end{equation}
for $1\le i,j\le n$. 

(ii) The map $g_{\bold x}: \T_b\to \T_b$ 
extends to a projective transformation of 
$P_-$ given by
\begin{equation}
\bold Z\to A_g(\bold x,b)\bold Z,
\end{equation} 
where  
$A_g$ is the matrix given by 
\begin{equation}
A_g(\bold x,b)_{ij}=X_{j-1}X_{i}^{-1}a, i\ge j;\  
A_g(\bold x,b)_{ij}=X_{j-1}X_{i}^{-1}b, i<j,  
\end{equation}
for $1\le i,j\le n$. 
\end{Proposition}

\begin{Proof} The proof is by a straightforward 
calculation. 
\end{Proof}

Now it is easy to prove Proposition \ref{res1}.
For this, we need to check the 
invertibility of the matrices $A_f$ and $A_g$. 
To prove the invertibility of this first matrix, it is sufficient to observe 
that $\text{det} A_f(\bold y,a)=(b-a)^{n-1}$ 
(this is easily established by induction).
The second matrix is the transpose of the first one, so the result follows. 

\section{Proof of Theorem \ref{Rmat} and Proposition \ref{res2}}

We start with proving Theorem \ref{Rmat}. 
Let us calculate the matrix
$A_f^{-1}(\bold y,a)$, computing the minors 
of $A_f$ of size $n-1$.  
We find that $(A_f^{-1})_{ij}=0$ unless $j-i$ is 0 or 1
modulo $n$, and 
\begin{equation}
(A_f^{-1})_{ii}=(b-a)^{-1}Y_iY_{i-1}^{-1},\
(A_f^{-1})_{i,i+1}=-(b-a)^{-1},\
(A_f^{-1})_{n,1}=-a(b-a)^{-1}.
\end{equation}
The matrix $A_g^{-1}$ is, as we mentioned, the transpose of $A_f^{-1}$. 

This leads to universal expressions for $f^{-1}$ and $g^{-1}$ which do not 
depend of $a,b$ and even of $n$, and have a ``local'' form
(in terms of the indices involved). These expressions will be useful below, 
but also seem interesting by themselves. 

Namely, let $\bold X,\bold Y$ denote infinite sequences 
consisting of all $X_i=T_i(\bold x),Y_i=T_i(\bold y)$, $i\in \Bbb Z$. 
Consider the functions from pairs of sequences to sequences, given by the 
formulas
\begin{equation}
\phi_{\bold Y}(\bold X)_i=Y_{i+1}Y_i^{-1}X_i-X_{i+1},
\gamma_{\bold X}(\bold Y)_i=(X_iX_{i-1}Y_i^{-1}-Y_{i-1}^{-1})^{-1}.
\end{equation}
The formulas for $A_f^{-1}$ and $A_g^{-1}$ imply the following.

\begin{Proposition} \label{univ} 
In terms of the projective coordinates $X_i,Y_i$, 
the maps $(\bold x,\bold y)\to f^{-1}_{\bold y}(\bold x)$ 
and $(\bold x,\bold y)\to g_{\bold x}^{-1}(\bold y)$ 
are given by the functions 
$\phi$ and $\gamma$ respectively.
\end{Proposition} 

Let us now use this proposition to show that $R$ is involutive.  

We must check that 
$\bold y=
g^{-1}_{f_{\bold y}^{-1}(\bold x)}\circ f^{-1}_{\bold x}(\bold y),
\ \bold x,\bold y\in 
(\C^*)^n$. 

According to Proposition \ref{univ}, 
this reduces to the identity 
\begin{equation}
(Y_{i+1}Y_i^{-1}X_i-X_{i+1})(Y_iY_{i-1}^{-1}X_{i-1}-X_i)^{-1}
(X_{i+1}X_i^{-1}Y_i-Y_{i+1})^{-1}-(X_iX_{i-1}^{-1}Y_{i-1}-Y_i)^{-1}=Y_i^{-1},
\end{equation}
which is straightforward. 

Now let is prove the quantum Yang-Baxter equation for $R$. 
It has three components. Let us prove the identity in the first components. 
It has the form
\begin{equation}
f_{\bold y}f_{\bold x}=f_{f_{\bold y}(\bold x)}f_{g_{\bold x}(\bold y)}.
\end{equation}
After changes of variables and inversions (using the involutivity of $R$) 
we see that this identity is equivalent to the claim that 
the map $f^{-1}_{f_{\bold y}^{-1}(\bold x)}f^{-1}_{\bold y}$
is a symmetric function of $\bold x,\bold y$. Applying this map to 
$\bold z$ using Proposition \ref{univ}, we find that this 
map is given by the operator on sequences defined by 
\begin{equation}
(B\bold Z)_i=(Y_{i+2}Y_{i+1}^{-1}X_{i+1}-X_{i+2})
(Y_{i+1}Y_i^{-1}X_i-X_{i+1})^{-1}
(Y_{i+1}Y_iZ_i-Z_{i+1})-(Y_{i+2}Y_{i+1}^{-1}Z_{i+1}-Z_{i+2}).
\end{equation}
But the symmetry of this operator in $\bold X,\bold Y$ is clear:
The coefficient of $Z_i$ is 
$\frac{Y_{i+2}X_{i+1}-Y_{i+1}X_{i+2}}{Y_{i+1}X_i-X_{i+1}Y_i}$, 
while the coefficient of $Z_{i+1}$ is
\begin{equation}
-(Y_{i+2}Y_{i+1}^{-1}X_{i+1}-X_{i+2})
(Y_{i+1}Y_i^{-1}X_i-X_{i+1})^{-1}-Y_{i+2}Y_{i+1}^{-1},
\end{equation}
which equals 
$-\frac{Y_{i+2}X_i-Y_iX_{i+2}}{Y_{i+1}X_i-X_{i+1}Y_i}$.
This proves the first component of QYBE.  

The third component of QYBE is completely analogous to the first one. 
The second component follows automatically 
from the involutivity and the first component, 
as shown in Proposition 2.2c in \cite{ESS} (the 
straightforward proof of this from \cite{ESS} applies verbatim 
in the case of rational maps). 
Theorem \ref{Rmat} is proved. 

Now let us prove Proposition \ref{res2}.
Recall that by definition, 
the group $G_R^+$ is a subgroup 
of $\text{Bir}(\T)$ generated by $f_{\bold y}^{-1}$, $\bold y\in \T$. 
In other words, $G_R^+$ is generated by the matrix functions 
$A_f(\bold z,\lambda)^{-1}$ of the parameter $\lambda\in \C^*$ for all 
$\bold z\in \T$.
Thus, $G_R^+\subset PGL_n(\C(\lambda))$. Similarly, 
$G_R^-\subset PGL_n(\C(\lambda))$ is generated by 
the matrix functions $A_g(\bold z,\lambda)$, $\bold z\in \T$.

Finally, it is easy to see that the groups $G_R^+,G_R^-$
are isomorphic, since $A_f(\bold z,\lambda)^T=A_g(\bold z,\lambda)$. 
 
Proposition \ref{res2} is proved. 

{\bf Remark.} It is shown in \cite{ESS} that for usual nondegenerate
involutive set-theoretical R-matrices, the groups 
$G_R^+,G_R^-$ are not only isomorphic but also conjugate 
in $\text{Aut}(X)$ by a certain bijection $T$, which intertwines
$f_x^{-1}$ and $g_x$. In our situation, however, the map 
$T$ is not well defined because of the ``rational'' 
character of the R-matrices. Therefore, the isomorphism 
of $G_R^+,G_R^-$ has to be checked independently.  

\section{Proof of Theorem \ref{comm}}

Consider an $n$ by $m$ matrix $(z_{kl})$ ($n,m>1$) with entries from $\C^*$.
Fix $i: 1\le i\le m-1$, and $j: 1\le j\le n-1$. 
We want to show that the involution $\sigma=PR_{j,j+1}$
acting on the rows of this matrix, commutes with the 
involution $\tau=PR_{i,i+1}$ acting on the columns. 
For this purpose, we will calculate the composition $\tau\sigma$, 
and show that it is ``symmetric under transposing the matrix'',
which would imply the commutativity. 

The elements which are subject to change under the composition $\tau\sigma$ 
are only those in rows $j,j+1$ or columns $i,i+1$. 
Therefore, it makes sense to introduce separate notations for these 
elements: $x_k=z_{jk},y_k=z_{j+1,k}$, $p_l=z_{li}, q_l=z_{l+1,i}$. 
Introduce $X_k,Y_k$ as in Section 4, and similarly $P_l,Q_l$. 
The same entries after application of $\sigma$, respectively 
$\tau\sigma$, will be denoted by the same letters with *, respectively 
**. Let $X_m=a,Y_m=b,P_n=c,Q_n=d$. The running subscript $k$ will be between
$0$ and $m$, and $l$ between $0$ and $n$. 
Since our rows and columns have 
four elements in common, we have the following four equations:
\begin{equation}
\label{proport}
\frac{P_j}{P_{j-1}}=\frac{X_i}{X_{i-1}};
\frac{Q_j}{Q_{j-1}}=\frac{X_{i+1}}{X_{i}};
\frac{Q_{j+1}}{Q_{j}}=\frac{Y_{i+1}}{Y_{i}};
\frac{P_{j+1}}{P_{j}}=\frac{Y_{i}}{Y_{i-1}}.
\end{equation}
Let us also introduce the notation 
\begin{equation}
G_k(\bold X,\bold Y)=
a\sum_{r=1}^k Y_r^{-1}X_{r-1}+b\sum_{r=k+1}^m Y_r^{-1}X_{r-1},
\end{equation}
and similarly
\begin{equation}
H_l(\bold P,\bold Q)=
c\sum_{r=1}^l Q_r^{-1}P_{r-1}+d\sum_{r=l+1}^n Q_r^{-1}P_{r-1}.
\end{equation}

According to section 4, we have 
\begin{equation}
X_k^*=X_k\frac{G_0(\bold X,\bold Y)}
{G_k(\bold X,\bold Y)}, 
\end{equation}
and 
\begin{equation}
Y_k^*=Y_k\frac{G_k(\bold X,\bold Y)}
{G_0(\bold X,\bold Y)},
\end{equation}
We also obviously have $P_l^*=P_l,Q_l^*=Q_l$ for $l\ne j$. 
Finally, for $l=j$, using formulas \ref{proport}, it is easy
to find that  
\begin{equation}\label{Peq}
P_j^*=P_j
\frac
{G_{i-1}(\bold X,\bold Y)}
{G_i(\bold X,\bold Y)}, 
\end{equation}
and 
\begin{equation}\label{Qeq}
Q_j^*=Q_j
\frac
{G_i(\bold X,\bold Y)}
{G_{i+1}(\bold X,\bold Y)}, 
\end{equation}
Now let us calculate the parameters after application
of $\tau\sigma$. We clearly have 
\begin{equation}
X_k^{**}=X_k^*=
X_k\frac{G_0(\bold X,\bold Y)}
{G_k(\bold X,\bold Y)}, k\ne i.
\end{equation}
and
\begin{equation}
Y_k^{**}=Y_k^*=
Y_k\frac{G_k(\bold X,\bold Y)}{G_0(\bold X,\bold Y)}, k\ne i
\end{equation}

Now let us compute $Q_l^{**}$. We have
\begin{equation}
Q_l^{**}=Q_l^*\frac{H_l(\bold P^*,\bold Q^*)}{H_0(\bold P^*,\bold Q^*)}.
\end{equation}
Let us express this via $P_r,Q_r$. As we know, $P_r^*=P_r$ and 
$Q_r^*=Q_r$, except for $r=j$. Thus, 
if $l\ne j$, all the stars in this expression 
can be removed, except in the combination 
$(Q_j^*)^{-1}P_{j-1}+Q_{j+1}^{-1}P_j^*$, which occurs both in the numerator 
and the denominator. 

Now, using the formulas for $P_j^*,Q_j^*$, we have
\begin{equation}
(Q_j^*)^{-1}P_{j-1}=Q_j^{-1}P_{j-1}+(a-b)
\frac{Q_j^{-1}P_{j-1}Y_{i+1}^{-1}X_i}
{G_i(\bold X,\bold Y)},
\label{stars1} 
\end{equation}
and
\begin{equation}
Q_{j+1}^{-1}P_j^*=Q_{j+1}^{-1}P_j-
(a-b)
\frac{Q_{j+1}^{-1}P_{j}Y_{i}^{-1}X_{i-1}}
{G_i(\bold X,\bold Y)}.
\label{stars2}
\end{equation}
But it is easy to check using equations \ref{proport} 
that 
\begin{equation}
Q_j^{-1}P_{j-1}Y_{i+1}^{-1}X_i= 
Q_{j+1}^{-1}P_{j}Y_{i}^{-1}X_{i-1}.\label{symme}
\end{equation}
This implies that 
\begin{equation}
(Q_j^*)^{-1}P_{j-1}+Q_{j+1}^{-1}P_j^*=
Q_j^{-1}P_{j-1}+Q_{j+1}^{-1}P_j,\label{nostars}
\end{equation}
and hence 
\begin{equation}
Q_l^{**}=Q_l\frac{H_l(\bold P,\bold Q)}{H_0(\bold P,\bold Q)}.
\end{equation}
We now notice the desired symmetry with the expression 
for $Y_k^{**}$, $k\ne i$, under the relabeling
$(P,Q,l,j,n,c,d,H)\to (X,Y,k,i,m,a,b,G)$. 
Similarly, 
\begin{equation}
P_l^{**}=
P_l\frac{H_0(\bold P,\bold Q)}{H_l(\bold P,\bold Q)},
\end{equation} 
which is symmetric to the expression for $X_k^{**}$. 

It remains to compute 
$X_i^{**},Y_i^{**}, P_j^{**}$ and $Q_j^{**}$ 
and make sure the desired symmetry 
is present. Using the symmetric images of the 
equations for \ref{Peq}, \ref{Qeq} for $P_j^*$, $Q_j^*$, we get 
\begin{equation}
X_i^{**}=X_i^*\frac
{H_{j-1}(\bold P^*,\bold Q^*)}
{H_j(\bold P^*,\bold Q^*)}=
\end{equation}
\begin{equation} 
X_i
\frac{G_0(\bold X,\bold Y)}
{G_i(\bold X,\bold Y)}
\cdot
\frac
{H_{j-1}(\bold P,\bold Q)}
{H_j(\bold P^*,\bold Q^*)}
\end{equation}
(the stars in the numerator are removed by formula \ref{nostars}). 
Now, by formulas \ref{stars1}, \ref{stars2}, we have 
\begin{equation}
H_j(\bold P^*,\bold Q^*)=
c\sum_{r=1}^{j} (Q^*_r)^{-1}P^*_{r-1}+d\sum_{r=j+1}^n (Q^*_r)^{-1}P^*_{r-1}=
\frac{S}{G_i(\bold X,\bold Y)},
\end{equation}
where 
\begin{equation}
S=(c-d)(a-b)U+G_i(\bold X,\bold Y)H_j(\bold P,\bold Q),
\ U=Q_j^{-1}P_{j-1}Y_{i+1}^{-1}X_i= 
Q_{j+1}^{-1}P_{j}Y_{i}^{-1}X_{i-1}. 
\end{equation}
Thus, 
\begin{equation}
X_i^{**}=X_iS^{-1}G_0(\bold X,\bold Y)H_{j-1}(\bold P,\bold Q).
\end{equation}
The element $P_j^{**}$ can now be found from equations \ref{proport}:
\begin{equation}
P_j^{**}=X_i^{**}P_{j-1}^{**}(X_{i-1}^{**})^{-1}=
P_jS^{-1}H_0(\bold P,\bold Q)G_{i-1}(\bold X,\bold Y).
\end{equation} 
Now it is apparent that $X_i^{**}$ is symmetric to 
$P_j^{**}$ under the symmetry
$(P,Q,l,j,n,c,d,H)\to (X,Y,k,i,m,a,b,G)$. 
Namely, it follows from the fact that the expression $U$, and hence 
$S$, is invariant 
under this symmetry. 

Finally, as we showed before, the expression
$X_i^{**}Y_i^{**}=X_{i-1}^{**}Y_{i-1}^{**}P_{j+1}^{**}/P_{j-1}^{**}$ 
is symmetric to 
$P_{j-1}^{**}Q_{j-1}^{**}X_{i+1}^{**}/X_{i-1}^{**}=P_j^{**}Q_j^{**}$, 
which implies that $Y_i^{**}$ is symmetric to $Q_j^{**}$. 
The theorem is proved.

\end{document}